\documentclass[12pt]{article}
\usepackage{amsmath}
\usepackage{amssymb}

\newtheorem{definition}{Definition}
\newtheorem{proposition}[definition]{Proposition}
\newtheorem{conjecture}[definition]{Conjecture}
\newtheorem{corollary}[definition]{Corollary}

\newcount\aaac
\def\aaa #1%
{%
\par
\medskip
\global\advance\aaac +1\relax
\noindent
\textbf{%
\the\aaac\relax
.\,\ %
#1%
}%
\hskip 0.5em\relax
\ignorespaces
}

\makeatletter
\@twosidetrue
\makeatother

\newcommand\DQ{\mathcal{D}_{\mathbb{Q}}}
\newcommand\Atilde{\tilde{\mathcal{A}}}
\newcommand\Aprime{\mathcal{A}'}
\newcommand\piG{\pi_{\mathcal{G}}}
\newcommand\piS{\pi_{S}}

\title{The algebra of balanced dessins}
\author{Jonathan Fine}
\date{}

\begin{document}

\maketitle

\aaa{} This paper gives a key definition, for a new approach to
dessins and algebraic numbers. The distant goal is to construct from
each dessin $D$ an algebraic number $\eta_D$, in a systematic and
useful way. The algebra of balanced dessins is generated by formal
sums $\psi_D$ of dessins, intended to be intermediate between $D$ and
$\eta_D$.

\aaa{} The wished for construction of $\eta_D$ has three phases. The
first is this. A \emph{dessin} is an ordered pair $(\alpha, \beta)$ of
permutations of a finite set $E$ (the edges), up to relabelling of the
edges. Cartesian product of edge sets induces a product on dessins.
We let $\mathcal{D}$ denote all irreducible dessins. We write
$1\in\mathcal{D}$ for the unique dessin with one edge.

\aaa{} We let $\DQ$ denote all rational coefficient finite formal sums
of elements of $\mathcal{D}$. The properties (i)~Cartesian product,
(ii)~decomposition into irreducible components, and (iii)~edge
relabelling, taken together induce a ring structure on $\DQ$.

\aaa{} For each $\psi\in\DQ$, the subalgebra generated by $\psi$ is
finite dimensional, as a vector space over $\mathbb{Q}$. (This follows
from biased dessins, see \cite{fine-dlbd}.) Thus, $\psi$ has a
\emph{minimal polynomial} $P_\psi(x)\in\mathbb{Q}[x]$. This is similar
to the minimal polynomial of an algebraic number, except that usually
$P_\psi(x)$ is reducible. The algebra $\DQ$ has divisors of zero.

\aaa{} The following shows that $\DQ$ is too large, for constructing
$\eta_D$.

\begin{proposition}
If $D\in\mathcal{D}$ then $P_D(x)$ is over $\mathbb{Q}$ a product of
linear factors.
\end{proposition}


\begin{corollary}
  Let $\mathcal{Z} \subseteq \DQ$ be an ideal. If $\mathcal{K} =
  \DQ/\mathcal{Z}$ is a field, then $\mathcal{K}=\mathbb{Q}$.
\end{corollary}

\aaa{} The second phase is to define $\psi_D\in\DQ$, having better
minimal polynomials. We use the symmetry group $\mathcal{G}$, the
\emph{absolute Galois group}, of the field $\mathcal{A}$ of all
algebraic numbers. It is already known that $\mathcal{G}$ acts on
$\mathcal{D}$, and thus on $\DQ$. This is a deep and central result
about dessins. Finding the $\mathcal{G}$-orbit of a dessin is usually
hard, and can be herculean. For more on this, see \cite[Chapter
  2]{LZ}.

\aaa{} We also use the symmetric group $S_3$. Each irreducible
dessin $D$ represents a finite cover $X_D\to\mathbb{P}_1(\mathbb{C})$
of the Riemann sphere, unramified away from $\{-1, +1,
\infty\}$. Permuting these three points induces an automorphism of the
sphere, and hence an action of $S_3$ on $\mathcal{D}$, and so on
$\DQ$. Finding the $S_3$-orbit is easy.

\aaa{} The key definition uses averages. For $\psi\in\DQ$ the orbit
$\mathcal{G}\psi$ is finite, because $\mathcal{G}$ preserves the
number of edges. Let $\piG(\psi)$ denote the average value (centre of
mass) of $\mathcal{G}\psi$, and similarly $\piS(\psi)$ for
$S_3\psi$. Both $\piG$ and $\piS$ are projections. They commute.

\begin{definition}
For $D\in\mathcal{D}$, we call $\psi_D= \piS(D)-\piS(\piG(D))$ a
\emph{balanced dessin}. We also say that $1\in\DQ$ is balanced. The
balanced dessins generate \emph{the algebra $\Atilde\subset\DQ$ of
  balanced dessins}.
\end{definition}

By design, $\psi_D$ is $S_3$-fixed, with $\mathcal{G}$-average
zero. If $D$ is $\mathcal{G}$-fixed, then $\psi_D=0$.

\aaa{} The third, final and hardest phase is this. To find an ideal
$\mathcal{Z}\subset\Atilde$ such that the quotient
$\Atilde/\mathcal{Z}=\Aprime$ is a useful field.  The best possible
outcome is $\Aprime \cong \mathcal{A}$.

\aaa{} Each element $\psi$ of $\Atilde$ has a minimal polynomial
$P(x)$. Thus, if $\Aprime$ is a field, then its elements are
algebraic numbers. We can say more.  Let $\eta\in\Aprime$ be the
residue of $\psi$. By definition $P(\psi)=0$, and so $P(\eta)=0$. When
$\Aprime$ is a field, we have $Z(\eta)=0$ for exactly one prime factor
$Z(x)$ of $P(x)$.

\aaa{} This leads to:

\begin{definition}
A \emph{choice of factors $Z$} is, for each $\psi \in \Atilde$, a
prime factor $Z_\psi(x)$ of the minimal polynomial $P_\psi(x)$.
We say that $Z$ is \emph{coherent} if $\Aprime = \Atilde/\mathcal{Z}$
is a field, where $\mathcal{Z}$ is the ideal $(Z_\psi(\psi) | \psi \in
\Atilde)$.
\end{definition}

\aaa{} Perhaps $\Atilde$ itself will guide us to our goal. The
condition $\Aprime\cong\mathcal{A}$ suggests:

\begin{conjecture}
For each $D\in\mathcal{D}$ and $\psi=\psi_D$, there is a unique prime
factor $Z_\psi(x)$ of $P_\psi(x)$ such that $\deg Z_\psi$ equals the
number of elements in $\mathcal{G}\psi$.
\end{conjecture}

\begin{conjecture}
This partial choice of factors has a unique coherent extension to
$\Atilde$.
\end{conjecture}

\aaa{} Should $\mathcal{Z}$ exist much as conjectured, then
$D\mapsto\psi_D\in\Atilde$ followed by $\Aprime = \Atilde/\mathcal{Z}$
should provide the desired construction of $\eta_D$ from $D$.

\aaa{} Such a $\mathcal{Z}$ would also raise questions such as (i)~are
$\Aprime$ and $\mathcal{A}$ isomorphic, (ii)~$Z_\psi(x)$ as a
$\mathcal{G}$-invariant of $D$, (iii)~construction via $\eta_D$ of the
Belyi pair $X_D\to \mathbb{P}_1(\mathbb{C})$, (iv)~the combinatorial
definition of $\mathcal{G}$, and (v)~efficient computation in
$\Aprime$.

\aaa{} Construction of $\mathcal{A}$ from $\mathcal{D}$ would also
perhaps open the door to showing that some special symmetries of
$\mathcal{D}$ (e.g.\@ \cite[\S3.10]{fine-dlbd}) also act on
$\mathcal{A}$, and so belong to $\mathcal{G}$.

\aaa{} For general background in dessins, see \cite[Chapter 2]{LZ} and
the references there. For more on the author's approach, see
\cite{fine-dlbd}.
I thank Joel Fine for helpful comments on an earlier draft of the
present paper.

\aaa{} The promise of useful coherence is such that, even if thought
unlikely or inaccessible, it is worth investigating. As usual,
$\Atilde$ is a profinite tower. The author intends to start by
calculating some basic examples \cite{py-dessins}.

\noindent Email: \texttt{jfine2358@gmail.com}\par
\noindent Location: Milton Keynes, UK\par
\noindent Date: 12 February 2018\par

\end{document}